\documentclass[11pt]{amsart}
\usepackage{tikz, tikz-3dplot}
\usepackage{amsmath, amsthm, amssymb}
\usepackage[letterpaper, top=3cm, bottom=3cm, left=3cm, right=3cm]{geometry}

\usepackage[T1]{fontenc}   
\usepackage[utf8]{inputenc} 

\usepackage{comment}
\usepackage{xcolor}
\usepackage{bm}  

\usepackage{hyperref}

\numberwithin{equation}{section}

\allowdisplaybreaks

\theoremstyle{plain}
\newtheorem{te}{Theorem}[section]

\newtheorem{lem}[te]{Lemma}
\newtheorem{co}[te]{Corollary}
\newtheorem{pr}[te]{Proposition}

\newtheorem{ex}[te]{Example}

\theoremstyle{remark}

\newtheorem*{ack*}{Acknowledgment}

\def\0{{\bf 0}}

\def\T{{\mathbb T}}
\def\R{{\mathbb R}}
\def\E{{\mathbb E}}

\def\C{{\mathbb C}}

\def\Z{{\mathbb Z}}
\def\P{{\mathbb P}}

\def\supp{{\operatorname{supp}\,}}
\def\nint{\mathop{\diagup\kern-13.0pt\int}}

\def\Pc{{\mathcal P}}

\begin{document}
	
\title{A remark on the majorant property for random sets}

\author{Hongki Jung}
	\address{Department of Mathematics, Louisiana State University, Baton Rouge, LA 70803}
\email{hjung@lsu.edu}

	
\begin{abstract}
We improve upon the majorant estimates introduced in \cite{MocSch} by applying Bourgain's $\Lambda(p)$ property to random sets.
		
\end{abstract}

\maketitle

\section{Introduction}
Let us write $e(x)=e^{2\pi i x}$ for $x\in \R$, and let $N\gg 1$ be a large integer. In \cite{HarLit}, Hardy and Littlewood showed that for an even integer $p$, one has
\begin{align}\label{eq:strict}
   \| \sum_{n\in A}a_n e(nx) \|_{L^p([0,1])} \leq \| \sum_{n\in A} e(nx) \|_{L^p([0,1])}
\end{align}
for any subset $A \subset [1,N]$,  $a_n\in \C$ and $|a_n|\leq 1$. We say that a subset $A\subset[1,N]$ has the (strict) majorant property if it satisfies \eqref{eq:strict}. 
When $p$ is not an even integer, subsets of $[1,N]$ that fail the majorant property by a factor of $N^\alpha$ for some $\alpha>0$ were discovered in \cite{GreRuz} and \cite{MocSch}. 

Nonetheless, a generic random subset still satisfies the majorant property. We construct such a random set as follows. Let $0<\delta<1$ and set $\tau=N^{-\delta}$. For $1\leq n\leq N$, let $\xi_n=\xi_n(\omega)$ be i.i.d. $\{0,1\}$-valued random variables with $\P(\xi_n=1)=\tau$, $\P(\xi_n=0)=1-\tau$. We define a random subset by
\begin{align*}
    S(\omega)=\{n\in[1,N]: \xi_n(\omega)=1\}.
\end{align*}
Mockenhaupt and Schlag proved the following majorant property in \cite{MocSch}. Here, $|a_n|$ denotes the absolute value of $a_n\in \{a_n\}_{n\in [1,N]}$.
\begin{te}
For any $\epsilon>0$ and $p\geq 2$, we have
    \begin{align*}
        \P \left( \sup_{|a_n|\leq 1} \|\sum_{n\in S(\omega)}a_ne(nx)\|_{L^p([0,1])}\geq N^\epsilon \|\sum_{n\in S(\omega)}e(nx)\|_{L^p([0,1])} \right) \rightarrow 0
    \end{align*}
as $N\rightarrow \infty$.
\end{te}

In this paper, we improve the estimate by removing the $N^\epsilon$ loss.
\begin{te}
    For any $0< \delta <1$, $p\geq 2$ and $\gamma>0$, there exist constants $C_{\gamma}$, $C_\gamma'$ such that
\begin{align*}
    \sup_{N\geq C_\gamma'} \; \P \left( \sup_{|a_n|\leq 1}\|\sum_{n\in S(\omega)}a_n e(nx)\|_{L^p([0,1])}\geq C_{\gamma}\|\sum_{n\in S(\omega)}e(nx) \|_{L^p([0,1])} \right) \leq \gamma.
\end{align*}
where the constants $C_\gamma$ and $C_\gamma'$ do not depend on $N$.
\end{te}

The case of random subsets $S(\omega)\subset [1,N]$ generated by stochastic processes was explored in \cite{TemSah}. In addition to subsets of $[1,N]$, certain subsets of the parabola and moment curves in $[1,N]^d$ that fail the strict majorant property \eqref{eq:strict} were investigated in \cite{BenBez} and \cite{GreGuo}. In Examples \ref{ex:randomsquare} and \ref{ex:randomparabola}, we follow a recent result from \cite{DeJuRy} to present subsets of $[1,N]^d$ whose random subsets satisfy the majorant property.

Furthermore, we improve the majorant estimate for a perturbation of arithmetic progressions introduced in \cite{MocSch}. For some $\epsilon_0 ,\epsilon_1>0$, let $L=N^{\epsilon_0}$, $s=N^{\epsilon_1}<a$ and define $I_j:=[(b+aj)-s, (b+aj)+s]$ for $1\leq j \leq L$. Let $\Pc := \cup_{j=1}^L I_j$. For $1\leq j \leq L$, let $\eta_j(\omega)$ be independent integer-valued random variables uniformly distributed on $I_j$. We define a random subset in $\Pc$, which represents a perturbation of arithmetic progressions, by
\begin{align*}
    S(\omega)=\{\eta_{j}(\omega)\}_{j=1}^L.
\end{align*}
In Section \ref{sectionlambda}, we show that the random subset $S(\omega)$ is a $\Lambda(p)$ set with high probability. As a consequence, we establish the majorant property for a perturbation of arithmetic progressions.

\bigskip

\subsection{Notations}
We write $A\lesssim B$ to denote $A\leq CB$ for positive quantities $A$ and $B$, and some constant $C>0$ throughout. If $A\lesssim B$ and $B\lesssim A$, then we write $A\sim B$. The indicator function of a set $E$ will be denoted by $\chi_E$ or $1_E$.

\begin{ack*}
    The author thanks Donggeun Ryou for insightful discussions.
\end{ack*}

\bigskip

\section{The Majorant Property for Random Sets}

For $N\gg 1$ and $0<\delta<1$, let $\tau=N^{-\delta}$. For $1\leq n\leq N$, let $\xi_n=\xi_n(\omega)$ be i.i.d. random variables with $\P(\xi_n=1)=\tau$, $\P(\xi_n=0)=1-\tau$. Here, $|a_n|$ denotes the absolute value of $a_n\in \{a_n\}_{n\in [1,N]}$. We define a random subset by
\begin{align}\label{eq:randomsubset}
    S(\omega)=\{n\in[1,N]: \xi_n(\omega)=1\}.
\end{align}

We prove that the majorant property holds with high probability for the random subsets $S(\omega)$ for all $0<\delta<1$ and all $p\geq2$ .
\begin{te}\label{te:main}
    For any $0< \delta <1$ and $p\geq 2$ we have
    \begin{align}\label{eqn:main}
        \E \left( \sup_{|a_n|\leq 1} \|\sum_{n\in S(\omega)}a_ne(nx)\|_{L^p([0,1])} \right) \leq C_p \E \left(\|\sum_{n\in S(\omega)}e(nx)\|_{L^p([0,1])} \right)
    \end{align}
for some $C_p$ not depending on $N$. Therefore, for any small $\gamma>0$, there exist constants $C_\gamma$, $C_\gamma'$ such that
\begin{align}\label{eqn:mainepsilon}
    \sup_{N\geq C_{\gamma}'} \; \P \left( \sup_{|a_n|\leq 1}\|\sum_{n\in S(\omega)}a_n e(nx)\|_{L^p([0,1])}\geq C_{\gamma}\|\sum_{n\in S(\omega)}e(nx) \|_{L^p([0,1])} \right) \leq \gamma.
\end{align}
\end{te}
We briefly show that \eqref{eqn:main} implies \eqref{eqn:mainepsilon}. By applying Chernoff-Hoeffding inequality (Lemma \ref{te:Lemma 2.2}) and the lower bound from Lemma \ref{te:Lemma 2.4}, we obtain
\begin{align*}
   \P \left (  \|\sum_{n\in S(\omega)}e(nx)\|_{L^p([0,1])} \leq \frac{1}{C_p} \E  (\|\sum_{n\in S(\omega)}e(nx)\|_{L^p([0,1])})  \right) \leq e^{-c\tau N}  
\end{align*}
for some $C_p>1$.
We apply Chebyshev's inequality to \eqref{eqn:main} to conclude that
\begin{align*}
&\P \left( \sup_{|a_n|\leq 1} \|\sum_{n\in S(\omega)}a_ne(nx)\|_{L^p([0,1])}\geq C_\gamma \|\sum_{n\in S(\omega)}e(nx)\|_{L^p([0,1])} \right)\\ &\leq
\P \left( \sup_{|a_n|\leq 1} \|\sum_{n\in S(\omega)}a_ne(nx)\|_{L^p([0,1])}\geq  \frac{C_\gamma}{C_p} \E (  \|\sum_{n\in S(\omega)}e(nx)\|_{L^p([0,1])} ) \right) +e^{-c\tau N}\\
&\leq \frac{C_p}{C_\gamma}+e^{-c\tau N} \leq \gamma
\end{align*}
where the last inequality holds provided that $C_\gamma$ and $\tau N$ are sufficiently large.

We recall a few preliminary results. The following lemma is classical, and its proof can be found in Lemma 4.1 of \cite{MocSch}.
\begin{lem}[Chernoff-Hoeffding]\label{te:Lemma 2.2}
 For fixed $0< \delta < 1$, let $\xi_n$ be random selectors and $\tau=N^{-\delta}$ as above. Then
 \begin{align*}
     \P\left( \sum_{n=1}^N\xi_n(\omega) \in [\frac{1}{2}\tau N ,2\tau N] \right)\geq 1- e^{-c\tau N}.
 \end{align*}
\end{lem}
Suppose that a sequence $\{a_{n}\}$ satisfies $|a_n|\leq 1$ for every $n\in [1,N]$. Observe that the lemma above implies $|a_n|_{l^2(n\in S(\omega))} \lesssim (\tau N)^{1/2}$ holds with high probability $\geq 1-e^{-c\tau N}$.

For $d\geq 1$, let $\Pc\subset [1,N]^{\alpha_1}\times [1,N]^{\alpha_2} \times \cdots \times [1,N]^{\alpha_d}\subset \Z^d$ with $|\Pc|=N$. We begin by proving a lower bound of the exponential sum over random subsets of $\Pc$.
\begin{lem}\label{te:Lemma 2.4}
 For fixed $0< \delta < 1$, let $\xi_n$ be random selectors of mean $\tau=N^{-\delta}$ as above. For $p\geq 2$, define
 \begin{align*}
     I_{p,N}(\omega)=\int_{[0,1]^d}|\sum_{n\in \Pc} \xi_n(\omega)e(n\cdot x)|^p dx.
 \end{align*}
 Then, we have
 \begin{align*}
     \E I_{p,N}(\omega) \gtrsim_p \tau^p N^{p-(\alpha_1+\cdots +\alpha_d)} +(\tau N)^{p/2}.
 \end{align*}
\end{lem}
\begin{proof}
    We apply Chernoff-Hoeffding's inequality followed by the constructive interference on $x \in [0,(100N)^{-\alpha_1}] \times [0, (100N)^{-\alpha_2}] \cdots \times [0,(100N)^{-\alpha_d}]$ to write
\begin{align*}
    \E \int_{[0,1]^d}|\sum_{n\in \Pc} \xi_n(\omega)e(n\cdot x)|^p dx
    &\geq \E \chi_{\sum_n \xi_n(\omega)\geq \frac{1}{2}\tau N}\int_{[0,1]^d}|\sum_{n\in \Pc} \xi_n(\omega)e(n\cdot x)|^p dx \\
    &\gtrsim (\tau N)^p N^{-(\alpha_1+\cdots \alpha_d)}.
\end{align*}

On the other hand, H\"older's inequality shows
\begin{align*}
    &\E \int_{[0,1]^d}|\sum_{n\in \Pc} \xi_n(\omega)e(n\cdot x)|^p dx\\
    &\geq \E \chi_{\sum_n \xi_n(\omega) \in [\frac{1}{2}\tau N, 2\tau N]}\int_{[0,1]^d}|\sum_{n \in \Pc} \xi_n(\omega)e(n \cdot x)|^p dx \\
    &\geq \E \chi_{\sum_n \xi_n(\omega) \in [\frac{1}{2}\tau N, 2\tau N]} (\int_{[0,1]^d}|\sum_{n\in \Pc} \xi_n(\omega)e(n \cdot x)|^2 dx)^{\frac{p}{2}} \\
    &\gtrsim  (\tau N)^\frac{p}{2}.
\end{align*}
    
\end{proof}

\begin{te}\label{te:Bourgain Lambda} For $0< \delta < 1$, let $\xi_n$ be random selectors of mean $\tau=N^{-\delta}$ as above. For $q=\frac{2(\alpha_1+\cdots \alpha_d)}{1-\delta}$, define
\begin{align*}
    K_q(\omega)=\sup_{|a_n|_{l^2_N}\leq 1}\|\sum_{n\in \Pc} a_n\xi_n(\omega)e(n \cdot x)\|_{L^q([0,1]^d)}.
\end{align*}
Assume that 
\begin{align}\label{eq:lambda(p)}
    \E K_q^q(\omega)\leq C_q
\end{align}
for some $0 <\delta <1$. Then, for such $\delta$, we have the random majorant property for any $p>2$
\begin{align*}
        \E \left( \sup_{|a_n|\leq 1} \|\sum_{n\in S(\omega)}a_ne(n\cdot x)\|_{L^p([0,1]^d)} \right) \leq C_p \E \left(\|\sum_{n\in S(\omega)}e(n\cdot x)\|_{L^p([0,1]^d)} \right)
    \end{align*}
\end{te}

\begin{proof}
For some $0< \delta <1$, write $q=\frac{2(\alpha_1+\cdots+\alpha_d)}{1-\delta}$ so that $\Lambda(q)$-property holds for generic random sets $S(\omega)$. We prove
\begin{align*}
        \E \left( \sup_{|a_n|\leq 1} \|\sum_{n\in S(\omega)}a_ne(n\cdot x)\|_{L^p([0,1]^d)}^p \right) \leq C_p \E \left(\|\sum_{n\in S(\omega)}e(n \cdot x)\|_{L^p([0,1]^d)}^p \right).
\end{align*}
Let us write $\E'$ for the restricted expectation
\begin{align*}
    \E' \sup_{|a_n|\leq 1} \|\sum_{n\in S(\omega)}a_ne(n \cdot x)\|_{L^p([0,1]^d)}^p :=\E\chi_{\sum \xi_n(\omega)\leq 2\tau N}(\omega) \sup_{|a_n|\leq 1} \|\sum_{n\in S(\omega)}a_ne(n \cdot x)\|_{L^p([0,1]^d)}^p 
\end{align*}
We use Lemma \ref{te:Lemma 2.2} to obtain
\begin{align*}
    \E  \sup_{|a_n|\leq 1} \|\sum_{n\in S(\omega)}a_ne(n\cdot x)\|_{L^p([0,1])}^p  &\leq  \E' \sup_{|a_n|\leq 1} \|\sum_{n\in S(\omega)}a_ne(n\cdot x)\|_{L^p([0,1]^d)}^p  + N^p e^{-c\tau N}\\
   &\leq \E' \sup_{|a_n|\leq 1} \|\sum_{n\in S(\omega)}a_ne(n \cdot x)\|_{L^p([0,1]^d)}^p  + O(1) .
\end{align*}
First, assume that $p\geq q$. We estimate the first term by
\begin{align*}
 &\E'\left( \sup_{|a_n|\leq 1} \int_{[0,1]^d} |\sum_{n=1}^N a_n 
 \xi_n(\omega) e(n\cdot x)|^p dx \right)\\   
 &\leq  \E' \left( \sup_{|a_n|\leq 1} \sup_{x\in[0,1]}|\sum_{n=1}^N a_n \xi_n(\omega) e(n \cdot x)|^{p-q}  \int_{[0,1]^d} |\sum_{n=1}^N a_n 
 \xi_n(\omega) e(n\cdot x)|^q dx \right)\\
 &\lesssim (\tau N)^{p-q} \E' \left( \sup_{|a_n|\leq 1}  \int_{[0,1]^d} |\sum_{n=1}^N a_n 
 \xi_n(\omega) e(n\cdot x)|^q dx \right).
\end{align*}
When estimating restricted expectation, note that $\sum \xi_n(\omega)\leq 2\tau N$, and hence $|a_n|_{l^2(n\in S(\omega))} \lesssim (\tau N)^{1/2}$. We further write
\begin{align*}
 \E'  \sup_{|a_n|\leq 1}  \int_{[0,1]^d} |\sum_{n=1}^N a_n 
 \xi_n(\omega) e(n \cdot x)|^q dx 
&\lesssim  \E'  \sup_{|a_n|_{l^2}\leq (\tau N)^{1/2}}  \int_{[0,1]^d} |\sum_{n=1}^N a_n 
 \xi_n(\omega) e(n\cdot x)|^q dx  \\
&\leq  \E  \sup_{|a_n|_{l^2}\leq (\tau N)^{1/2}}  \int_{[0,1]^d} |\sum_{n=1}^N a_n 
 \xi_n(\omega) e(n \cdot x)|^q dx  \\
& \lesssim (\tau N)^{q/2}
\end{align*}
where the last inequality follows from \eqref{eq:lambda(p)}. Combining the last two inequalities, we obtain
\begin{align*}
(\tau N)^{p-q} \E' \left( \sup_{|a_n|\leq 1}  \int_{[0,1]^d} |\sum_{n=1}^N a_n 
 \xi_n(\omega) e(n\cdot x)|^q dx \right) &\lesssim (\tau N)^{p-q} (\tau N)^{q/2} = (\tau N)^{p-\frac{q}{2}}.
\end{align*}
Using the fact $(\tau N)^{-q/2}=N^{-(\alpha_1+\cdots \alpha_d)}$ and Lemma \ref{te:Lemma 2.4}, we conclude that
 \begin{align*}
     (\tau N)^{p-\frac{q}{2}} = (\tau N)^p N^{-(\alpha_1+\cdots \alpha_d)} \lesssim \E I_{p,N}(\omega) = \E \left(\|\sum_{n\in S(\omega)}e(n\cdot x)\|_{L^p([0,1]^d)}^p \right).
 \end{align*}

For the case $p < q$, we apply H\"older's inequality to write
\begin{align*}
    \E'\left( \sup_{|a_n|\leq 1} \int_{[0,1]^d} |\sum_{n=1}^N a_n 
 \xi_n(\omega) e(n\cdot x)|^p dx \right) &\leq \E'\left( \sup_{|a_n|\leq 1} \int_{[0,1]^d} |\sum_{n=1}^N a_n 
 \xi_n(\omega) e(n\cdot x)|^q dx \right)^{p/q}\\
  &\leq \E\left( \sup_{|a_n|_{l^2}\leq (\tau N)^{1/2}} \int_{[0,1]^d} |\sum_{n=1}^N a_n 
 \xi_n(\omega) e(n\cdot x)|^q dx \right)^{p/q}.
\end{align*}
Apply \eqref{eq:lambda(p)} and Lemma \ref{te:Lemma 2.4} consecutively to conclude that
\begin{align*}
    \E\left( \sup_{|a_n|_{l^2}\leq (\tau N)^{1/2}} \int_{[0,1]^d} |\sum_{n=1}^N a_n 
 \xi_n(\omega) e(n\cdot x)|^q dx \right)^{p/q}
 &\lesssim (\tau N)^{p/2} \lesssim \E I_{p,N}(\omega)\\
 &= \E \|\sum_{n\in S(\omega)}e(n\cdot x)\|_{L^p([0,1]^d)}^p .
\end{align*}

\end{proof}

\begin{ex}\label{ex:randomsubset}
The random subsets $S(\omega)\subset [1,N]$, defined by \eqref{eq:randomsubset}, satisfy condition \eqref{eq:lambda(p)} for $0<\delta<1$, as shown in \cite{Bo}. Therefore, the majorant property holds for all $0<\delta<1$ and $p>2$. This completes the proof of Theorem \ref{te:main}.
\end{ex}

\begin{ex}\label{ex:randomsquare} The random subsets of squares $S(\omega)=\{  n^2 : 1\leq n \leq N, \; \xi_n(\omega)=1 \}\subset [1,N^2]$ satisfy condition \eqref{eq:lambda(p)} for $0<\delta<1$ as shown in Section 4 of \cite{Bo2} and Theorem 5.2 of \cite{DeJuRy}.  Thus, the random majorant property holds for all $0<\delta<1$ and $p>4$.
\end{ex}

\begin{ex}\label{ex:randomparabola}
In higher dimensions, random subsets of parabola $\{(n,n^2)\}\subset [1,N]\times [1,N^2]$ and paraboloid $\{(n,m,n^2+m^2)\}\subset [1,N^{1/2}]^2\times [1,N]$ satisfy condition \eqref{eq:lambda(p)} by following the approach in \cite{DeJuRy}, which combines equation (3.4) from \cite{Bo3} and Theorem 5.2 from \cite{DeJuRy}. Consequently, the random majorant property holds for any $0< \delta <1$ and $p>6$ and $p>4$, respectively.
\end{ex}

The following corollary improves Theorem 4.16 from \cite{MocSch} by applying the improved estimate in Example \ref{ex:randomsubset}. Let $\omega \in \T=\R/\Z$ and let the probability measure $\P(d\omega)$ be the Lebesgue measure on $[0,1]$. 
We define correlated selectors $\xi_j(\omega):=1_{[0,\tau]}(2^j \omega)$ for $j\leq 1$, where $0< \tau < 1$.

\begin{co}
    Let $N\gg 1$ be large integer, $0< \delta < 1$ be fixed and $\tau=N^{-\delta}$.  Let $\xi_j$ be correlated selectors as above and we define a subset
    \begin{align*}
        S(\omega)=\{j\in  [1,N]: \xi_j(\omega)=1\}
    \end{align*}
    for every $\omega \in \T$. Then we have
        \begin{align*}
        \E  \sup_{|a_n|\leq 1} \|\sum_{n\in S(\omega)}a_ne(nx)\|_{L^p([0,1])}^p  \lesssim_p (\log N) \; \E \|\sum_{n\in S(\omega)}e(nx)\|_{L^p([0,1])}^p .
    \end{align*}
\end{co}

\bigskip

\section{The Majorant Property for a Perturbation of Arithmetic Progressions}

In this section, we establish the majorant property for a perturbation of arithmetic progressions introduced in Section 5 of \cite{MocSch}. For some $\epsilon_0, \epsilon_1>0$, let $L=N^{\epsilon_0}$, $s=N^{\epsilon_1}<a/2$ and define an arithmetic progression of length $L$ by
\begin{align*}
  \{b+aj :  1 \leq j \leq L \}\subset [1,N].
\end{align*}
for some $0< b < a$.
Define the intervals $I_j:=[(b+aj)-s, (b+aj)+s]$ for $1\leq j \leq L$ and their union by $\Pc := \bigcup_{j=1}^L I_j$. For $1\leq j \leq L$, let $\eta_j(\omega)$ be an independent integer-valued random variables uniformly distributed on $I_j$. We define a random subset of $\Pc$ by
\begin{align}\label{eq:spap}
    S(\omega)=\{\eta_{j}(\omega)\}_{j=1}^L.
\end{align}

The corresponding $\Lambda(p)$ property for this random set will be proved in Section \ref{sectionlambda}. Note that for $p={2(\epsilon_0+\epsilon_1)/\epsilon_0}$, we have $L=(Ls)^{2/p}$, which is the maximal size of $\Lambda(p)$ subsets of $\Pc$. 
\begin{te}[$\Lambda(p)$ Property]\label{te:lambdapap}
    Let $S(\omega)$ and $\eta_j(\omega)$ be defined as above. For $p=\frac{2(\epsilon_0+\epsilon_1)}{\epsilon_0}$, define
\begin{align*}
    K_p(\omega)=\sup_{|a_n|_{l^2_N}\leq 1}\|\sum_{n\in S(\omega)} a_ne(nx)\|_{L^p([0,1])}.
\end{align*}
Then $\E K_p^p(\omega)\leq C_p$.
\end{te}

To prove the majorant property, we first estabilish a lower bound for the exponential sum over random sets.
\begin{lem}\label{te:lowerboundpap}
 Let $q>2$ be defined by equation $s=L^{q/2-1}$. For $p> 2$, define
\begin{align*}
     I_{p,N}(\omega)=\int_{[0,1]}|\sum_{n\in S(\omega)} e(n x)|^p dx.
 \end{align*}
 Then, we have
 \begin{align*}
    \E I_{p,N}(\omega) \gtrsim_p  \frac{L^{p-1}}{s} +L^{\frac{p}{2}}  .
 \end{align*}
\end{lem}
\begin{proof}
First, H\"older's inequality shows that for any $\omega\in \Omega$,
\begin{align*}
    \int_{[0,1]}|\sum_{n\in S(\omega)} e(n x)|^p dx \geq  (\int_{[0,1]}|\sum_{n\in S(\omega)} e(n  x)|^2 dx)^{\frac{p}{2}} =  L^\frac{p}{2}.
\end{align*}

If $\frac{L^{p-1}}{s} > L^{p/2}$, we have $L^{p/2-1} > s = L^{q/2-1}$. Choose $r$ such that $q< r< p$. We write

\begin{align*}
    & I_{p,N}(\omega) \\
    &\geq \int_0^{1/s} |\sum_{n\in \Z}\E1_{S(\omega)}(n)e(nx)|^p dx- \int_0^{1/s} |(\sum_{n\in \Z} 1_{S(\omega)}(n) - \E1_{S(\omega)}(n))e(nx)|^p dx \\
    &=: I_1 - I_2.
\end{align*}

Then $I_1\gtrsim L^{p-1}/s$ holds by equation (5.18) in \cite{MocSch}. In order to estimate $I_2$, we set $\lambda:=L^{\frac{r-q}{2p}}$ and apply Lemma 5.3 and Corollary 5.4 from \cite{MocSch} to write
\begin{align*}
   \P \left( \|(\sum_{n\in \Z} 1_{S(\omega)}(n) - \E1_{S(\omega)}(n))e(nx) \|_{L^\infty_x} > \lambda \sqrt{L} \right) < C(s+L)^2e^{-c\lambda^2}.
\end{align*}
Therefore, we have
\begin{align*}
   \E I_2  \lesssim_p \frac{1}{s}\lambda^pL^{\frac{p}{2}} + \frac{1}{s}L^p (s+L)^2e^{-c\lambda^2}\lesssim \frac{1}{s} L^{\frac{p+r-q}{2}} \ll \frac{1}{s}L^{p-1}.
\end{align*}
Combining the estimates for $I_1$ and $\E_\omega I_2$, we get
\begin{align*}
    \E I_{p,N} \gtrsim_p \frac{L^{p-1}}{s}+ L^{\frac{p}{2}}.
\end{align*}

\end{proof}

Following the previous approach, we prove the majorant property for the random set \eqref{eq:spap}.

\begin{te}
    For $S(\omega)$ defined in \eqref{eq:spap}, we have
    \begin{align*}
        \E \left( \sup_{|a_n|\leq 1} \|\sum_{n\in S(\omega)}a_ne(nx)\|_{L^p([0,1])} \right)
        \leq C_p \E \left(\|\sum_{n\in S(\omega)}e(nx)\|_{L^p([0,1])} \right)
    \end{align*}
for some $C_p$ not depending on $N$. Therefore, for any $\gamma>0$, there exist constants $C_\gamma$, $C_\gamma'$ such that
\begin{align*}
    \sup_{N\geq C_\gamma'} \; \P \left( \sup_{|a_n|\leq 1}\|\sum_{n\in S(\omega)}a_n e(nx)\|_{L^p([0,1])}\geq C_{\gamma}\|\sum_{n\in S(\omega)}e(nx) \|_{L^p([0,1])} \right) \leq \gamma.
\end{align*}
\end{te}

\begin{proof}
Let $q>2$ be defined by equation $s=L^{q/2-1}$, so that $L=(Ls)^{2/q}$, and $\Lambda(q)$-property holds for generic random subsets $S(\omega) \subset \Pc$. First, assume that $p> q$. We compute
\begin{align*}
 &\E \sup_{|a_n|\leq 1} \int_{[0,1]} |\sum_{n\in S(\omega)} a_n 
  e(nx)|^p dx \\
 &\leq  \E  \sup_{|a_n|\leq 1} (\sup_{x\in[0,1]}|\sum_{n\in S(\omega)} a_n  e(nx)|^{p-q} ) \int_{[0,1]} |\sum_{n\in S(\omega)} a_n 
 e(nx)|^q dx \\
 &\lesssim L^{p-q} \E \left( \sup_{|a_n|\leq 1}  \int_{[0,1]} |\sum_{n\in S(\omega)} a_n  e(nx)|^q dx \right).
\end{align*}
Note that $|a_n|_{l^2(n\in S(\omega))} \leq L^{1/2}$. We further write
\begin{align*}
 \E \left( \sup_{|a_n|\leq 1}  \int_{[0,1]} |\sum_{n\in S(\omega)} a_n  e(nx)|^q dx \right) 
&\lesssim  \E \left( \sup_{|a_n|_{l^2}\leq L^{1/2}}  \int_{[0,1]} |\sum_{n\in S(\omega)} a_n 
  e(nx)|^q dx \right) \\
& \lesssim L^{q/2}
\end{align*}
where the last inequality follows from Theorem \ref{te:lambdapap}. Combining the last two inequalities, we obtain
\begin{align*}
L^{p-q} \E \left( \sup_{|a_n|\leq 1}  \int_{[0,1]} |\sum_{n\in S(\omega)} a_n 
  e(nx)|^q dx \right) &\lesssim L^{p-\frac{q}{2}}= \frac{L^{p-1}}{s} \lesssim \E I_{p,N}(\omega).
\end{align*}
using $L^{-q/2}=(Ls)^{-1}$ and Lemma \ref{te:lowerboundpap}.

For the case $p \leq q$, we apply H\"older's inequality to write
\begin{align*}
    \E\left( \sup_{|a_n|\leq 1} \int_{[0,1]} |\sum_{n\in S(\omega)} a_n e(nx)|^p dx \right) &\leq \E \left( \sup_{|a_n|\leq 1} \int_{[0,1]} |\sum_{n\in S(\omega)} a_n e(nx)|^q dx \right)^{p/q}\\
  &\leq \E\left( \sup_{|a_n|_{l^2}\leq L^{1/2}} \int_{[0,1]} |\sum_{n\in S(\omega)} a_n e(nx)|^q dx \right)^{p/q}.
\end{align*}
Apply Theorem \ref{te:lambdapap} and Lemma \ref{te:lowerboundpap} consecutively to conclude that
\begin{align*}
    \E\left( \sup_{|a_n|_{l^2}\leq L^{1/2}} \int_{[0,1]} |\sum_{n\in S(\omega)} a_n e(nx)|^q dx \right)^{p/q}
 &\lesssim L^{p/2} \lesssim  \E I_{p,N}(\omega).
\end{align*}

\end{proof}

\bigskip

\section
{ The \texorpdfstring{$\Lambda(p)$}{}  Property for a Perturbation of Arithmetic Progressions}\label{sectionlambda}

We establish the $\Lambda(p)$ property in a more general setting. Let $\Phi=(\varphi_i)_{i=1}^N$ be a sequence of mutually orthogonal real-valued functions satisfying the uniform bound $\|\varphi_i\|_{L^\infty}\leq 1$ for each $i\in [1,N]$. For a given $p>2$, let $L=N^{2/p}$, $s=N/L$ and $\{I_j\}_{j=1}^L$ be the partition of $[1,N]$ with $|I_j|=s$. Let $\{\eta_j\}_{j=1}^L$ be independent, integer-valued random variables uniformly distributed over $I_j$. Define the random subset 
\begin{align}\label{eq:slambda}
S(\omega)=S_\omega=\{\eta_j(\omega)\}_{j=1}^L.
\end{align}
Observe that the random subset $S_\omega$ selects one element from each block $I_j$, $1\leq j \leq L$.

For $S\subset [1,N]$ define  
\begin{equation}\label{eq:KS}
K(S) := \sup_{|\mathbf{a}|_{l^2}\le 1}\|\sum_{i\in S}a_i\varphi_i\|_{L^p(dx)},
\end{equation}
where $\mathbf a = (a_1,\dots,a_N)$.  Let $K(\omega) = K(S_\omega)$. We follow the detailed exposition in \cite{JuLaOrVu} to prove that there exists $C(p)>0$ depending only on $p$ such that 
\begin{equation}\label{eq:random-lambdap}
    \mathbb E K(\omega)^p \le C(p).
\end{equation}

\bigskip

\subsection{The case \texorpdfstring{$2<p<4$}{}}
For given $p>2$, let $L=N^{2/p}$, $s=N/L$ and let $\{I_j\}_{j=1}^L$ be a partition of $[1,N]$ with $|I_j|=s$.
Let $\eta_j: \Omega \rightarrow I_j$ be independent random variables on some probability space $\Omega$ satisfying the uniform distribution
\begin{align*}
\P(\eta_j=i)=\frac{1}{s} \quad \text{for} \quad i\in I_j.
\end{align*}
Consider the random set
\[
S_\omega= S(\omega) = \{\eta_j(\omega): j\in [1,L]\}, \qquad \omega\in\Omega,
\]
which has cardinality $|S_\omega|= L =N^{2/p}$. For $\omega\in\Omega$ and a sequence $\mathbf{d}=(d_1,\dots,d_N)$, define
\begin{align*}
f_{\mathbf{d},\omega}=\sum_{i\in S(\omega)}d_i\varphi_i.
\end{align*}

We substitute $n=N$, $n_0=N^{2/p}=L$ into Proposition 3.11 from \cite{JuLaOrVu} to obtain the following proposition.
\begin{pr}\label{prop:dec}
If $p>2$, then
\begin{align*}
\int_\Omega K^p(\omega)\,d\omega&\lesssim 
1+\int_\Omega K^{p-1}(\omega)\,d\omega
\\
&\quad+\int_\Omega\int_\Omega\int_\Omega \sup_{(\bm a,\bm b,\bm c)\in \mathcal{A}_1}\bigg|\langle f_{\bm a,\omega_1}, 	f_{\bm b,\omega_2}(1+|f_{\bm c,\omega_3}|)^{p-2}\rangle\bigg|\,d\omega_1\,d\omega_2\,d\omega_3
\\
&\quad+\int_\Omega\int_\Omega\int_\Omega \sup_{(\bm a,\bm b,\bm c)\in \mathcal{A}_2}\bigg|\langle f_{\bm a,\omega_1}, 	f_{\bm b,\omega_2}(1+|f_{\bm c,\omega_3}|)^{p-2}\rangle\bigg|\,d\omega_1\,d\omega_2\,d\omega_3,
\end{align*}
where the suprema are taken over the sets 
\begin{align*}
\mathcal{A}_1:=\{(\bm a,\bm b,\bm c):\ &|\supp \bm a|, |\supp \bm b|, |\supp \bm c|\le n_0, |\bm a|,|\bm b|,|\bm c|\le 1
\\ 
&\max_{1\le i\le N}|a_i|\le (|\supp\bm b|+|\supp\bm c|)^{-1/2}
\\
&\quad\text{and}\quad \max_{1\leq j \leq L} (|\supp \bm a\cap I_j|, |\supp \bm b\cap I_j|, |\supp \bm c\cap I_j|) \leq 1\},
\\
\mathcal{A}_2:=\{(\bm a,\bm b,\bm c):\ &|\supp \bm a|, |\supp \bm b|, |\supp \bm c|\le n_0, |\bm a|,|\bm b|,|\bm c|\le 1 ,
\\
&\max_{1\le i\le N}(|a_i|,|b_i|)\le |\supp\bm c|^{-1/2}
\\
&\quad \text{and}\quad \max_{1\leq j \leq L} (|\supp \bm a\cap I_j|, |\supp \bm b\cap I_j|, |\supp \bm c\cap I_j|) \leq 1 \}.
\end{align*}
\end{pr}

\begin{proof}  Let $\{\eta_i\}_{i=1}^N,\{\zeta_i\}_{i=1}^N$ be independent $\{0,1\}$-valued random variables of respective means
\begin{equation*}
\int\eta_i(t)\,dt = \frac13 \quad\text{and}\quad \int\zeta_i(t)\,dt = \frac12, \qquad 1\le i\le N.
\end{equation*}
Define the disjoint random sets
\begin{align*}
R^1_t &= \{1\le i \le N: \eta_i(t) = 1\},\\
R^2_t &= \{1\le i\le N : \eta_i(t) = 0, \zeta_i(t) = 1\},\\ 
R^3_t &= \{1\le i\le N : \eta_i(t) = 0, \zeta_i(t) = 0\}.
\end{align*}
Observe that for any $\omega\in \Omega$, $S=S_\omega$ satisfies $\max_{1\leq j \leq L}|S\cap I_j|=1$. Accordingly, we replace definition $\mathcal{A}(S)$ in Proposition 3.6 from \cite{JuLaOrVu} by

\begin{align*}
\mathcal{A}(S):=\{(\mathbf{a},\mathbf{b},I):& I\subset [1,N], \mathbf{a} = (a_i)_{i\in I\cap S} , \mathbf{b}=(b_i)_{i\in S}, |\mathbf{a}|,|\mathbf{b}|\le 1, \max_i|b_i| \le |I|^{-1/2}\\
&\quad \text{and}\quad \max_{1\leq j \leq L} (|\supp \bm a\cap I_j|, |\supp \bm b\cap I_j|, |\supp \bm c\cap I_j|) \leq 1 \}.
\end{align*}
Then, we substitute the expressions of $J_1$ and $J_2$  by
\begin{align*}
J_1:= &\quad\iint \sup_{(\mathbf{a},\mathbf{b},I)\in \mathcal{A}(S_\omega)} |\langle \sum_{i\in R_t^1\cap S(\omega)}b_i\varphi_i, (\sum_{i\in I\cap R_t^2\cap S(\omega)}a_i\varphi_i) 
\\
&\qquad\qquad\qquad\qquad\qquad\times(1+ |\sum_{i\in I\cap R_t^3\cap S(\omega)}a_i\varphi_i|)^{p-2}\rangle|\,d\omega\,dt
\\
J_2:=&\quad\iint\sup_{(\mathbf{a},\mathbf{b},I)\in \mathcal{A}(S_\omega)}|\langle\sum_{ i\in R^1_t\cap S(\omega)} b_i\varphi_i,(\sum_{i \in R^2_t\cap S(\omega)}b_i\varphi_i)
\\
&\qquad\qquad\qquad\qquad\qquad\times(1+|\sum_{i \in I\cap R^3_t\cap S(\omega)}a_i\varphi_i|)^{p-2}\rangle| \,d\omega\,dt.
\end{align*}
and apply the proof of Proposition 3.11 from \cite{JuLaOrVu}.

\end{proof}

Let $q_0 = \log N$ and for $1\le m\le N$ define 
\begin{align}\label{eq:Pim}
&		\Pi_{m}=\{\bm a=(a_i)_{i=1}^N:|\bm a|\le 1, |\supp \bm a|\le m \quad\text{and} \quad \max_{1\leq j\leq N}|\supp \bm a \cap I_j|\leq 1 \}.		
\end{align}
For $\omega_1, \omega_2$, $\omega_3\in\Omega$ and $m_1,m_2,m_3\in[1,N]$, define
\begin{align}\label{eq:K_{m_1,m_2,m_3}}
&	K_{m_1,m_2,m_3}(\omega_1,\omega_2,\omega_3)=\sup_{\substack{|A|\le m_1\\\max_{1\leq j\leq L}|A\cap I_j|= 1}}\sup_{\bm b\in \Pi_{m_2}}\sup_{\bm c\in \Pi_{m_3}}\frac{1}{\sqrt{m_1}}\sum_{i\in A\cap S(\omega_1)} |\langle \varphi_i, 	f_{\bm b,\omega_2}(1+|f_{\bm c,\omega_3}|)^{p-2}\rangle |.
\end{align}
For $2<p<4$, the following estimate holds for any $\omega_2$ and $\omega_3$
\begin{align*}
	\|K_{m_1,m_2,m_3}(\omega_1,\omega_2,\omega_3)\|_{L^{q_0}(d\omega_1)}\le (\delta m_3^{\frac{p}{2}-1}+\frac{m_2+m_3}{m_1})^{\frac{1}{2}}(1+K(\omega_2)+K(\omega_3))^{p-\sigma}
\end{align*}
where $ \delta=N^{\frac{2}{p}-1}$, as stated in Theorem \ref{te:k2p4}.

We substitute equation (3.21) in \cite{JuLaOrVu} for $I_1$, $I_2$ adapted to $\mathcal{A}_1$, $\mathcal{A}_2$ arising in Proposition \ref{prop:dec} by
\begin{align} \label{eq:I12def}
    I_j &= \iiint \sup_{(\bm a,\bm b, \bm c)\in \mathcal A_j}\left|\langle f_{\bm a,\omega_1}, f_{\bm b,\omega_2}(1+|f_{\bm c,\omega_3}|)^{p-2}  \rangle \right|\,d\omega_1\,d\omega_2\,d\omega_3, \qquad j=1,2.
\end{align}
Then, we apply the proof in Section 3.4 of \cite{JuLaOrVu} showing that  $\|K(\omega)\|_{L^p(d\omega)}\lesssim 1$.  

\bigskip

\subsection{The case \texorpdfstring{$p>4$}{}.}\label{subsec:p>4}

We apply Bourgain's inductive argument, following the exposition in Section 3.5 of \cite{JuLaOrVu}. We aim to prove that, for a given partition $\{\bar{I_j}\}_{j=1}^{N^{2/p}}$ of $[1,N]$ with $|\bar{I}_j|=N^{1-2/p}$,  a random subset $S(\omega)\subset[1,N]$ of size $L=N^{2/p}$, as defined in \eqref{eq:slambda}, satisfies the $\Lambda(p)$-property. 

Assume that \eqref{eq:random-lambdap} holds for some exponent $p_1$ such that $p_1 < p < 2p_1$. Let $L_1=N^{2/p_1}$ and $\{I_k\}_{k=1}^{L_1}$ be a refinement of $\{\bar{I_j}\}_{j=1}^{L}$ with $|I_k|=N^{1-2/p_1}$. Let $\eta_k: \Omega \rightarrow I_k$ be independent random variables on a probability space $\Omega$ satisfying the uniform distribution
\begin{align*}
\P(\eta_k=i)=N^{\frac{2}{p_1}-1} \quad \text{for} \quad i\in I_k  .  
\end{align*}
By the induction hypothesis, we assume that $S_1=\{\eta_k(\bar{\omega})\}_{k=1}^{L_1}$ is a random subset satisfying the $\Lambda(p_1)$-property
\begin{equation}\label{eq:lambda-p1}
  \sup_{|\bm{a}|_{l^2}\leq 1} \|\sum_{i\in S_1} a_i\varphi_i\|_{p_1}\le C.
\end{equation}

Let $s'=L_1/L= N^{2/p_1-2/p}$. For each $1\leq j \leq L$, let $\{k(j)\}\subset\{k : 1 \leq k \leq L'\}$ be a subset of indices such that $I_{k(j)}\subset \bar{I_j}$ and $|\{k(j)\}|=s'$. Let $\zeta_j: \Omega' \rightarrow \{k(j)\}$ be independent random variables on a probability space $\Omega'$, uniformly distributed by
\begin{align*}
\P(\zeta_j=k(j))=\frac{1}{s'}.  
\end{align*}
Let $\{I_j'\}_{j=1}^{L}$ be a partition of $S_1$, where $I_j'=\{\eta_{k(j)}(\bar{\omega})\}_{k(j)}$ and $|I_j'|=s'$.

Now, define a random subset $S_\omega\subset S_1$ using the random variables $\{\zeta_j\}_{j=1}^L$ by
\begin{align}\label{eq:sp4}
    S_\omega=\{\eta_{\zeta_j(\omega)} (\bar{\omega})  :  1\leq j \leq L \}.
\end{align}
We define $K_{S_1}(\omega)$ by restricting the summation in the definition of $K_{S_1}$ to $S_\omega$, specifically:
$$
K_{S_1}(\omega) = \sup_{|\mathbf a|_{l^2}\le 1}\|\sum_{i\in S_\omega}a_i\varphi_i\|_{p}.
$$

Assuming that \eqref{eq:lambda-p1} holds for $S_1\subset[1,N]$ of size $L_1$ with $p_1$ satisfying $p/2< p_1 < p$, we will prove in Theorem \ref{te:kp4} the following analog of Theorem \ref{te:k2p4}
\begin{equation}\label{eq:Kest-p>4}
\|K^{S_1}_{m_1,m_2,m_3}(\omega_1,\omega_2,\omega_3)\|_{L^{q_0}(d\omega_1)}\lesssim (\delta' m_3^{p/p_1-1}+\frac{m_2+m_3}{m_1})^{1/2}(1+K_{S_1}(\omega_2)+K_{S_1}(\omega_3))^{p-\sigma},
\end{equation}
with $\delta'=N^{2/p-2/p_1}$, where 
$$
K^{S_1}_{m_1,m_2,m_3}(\omega_1,\omega_2,\omega_3):=\sup_{\substack{A\subset S_1 \\ |A|\le m_1\\ \max_{1\leq j \leq L} |A\cap I_j'|\leq 1}} \sup_{\bm b\in \Pi^{S_1}_{m_2}}\sup_{\bm c\in \Pi^{S_1}_{m_3}}\frac{1}{\sqrt{m_1}}\sum_{i\in A\cap S(\omega_1)} |\langle \varphi_i, 	f_{\bm b,\omega_2}|f_{\bm c,\omega_3}|^{p-2}\rangle |
$$
and
\begin{align*}
&\Pi^{S_1}_{m}=\{\mathbf{a}=(a_i)_{i\in S_1}:|\mathbf{a}|\le 1, |\supp \mathbf{a}|\le m\ \quad \text{and}\quad \max_{1\leq j \leq L} |\supp \bm a \cap I_j'|\leq 1\}.		
\end{align*}

We substitute $n=N$ and $n_0=N^{2/p}$, and adjust the vectors to be supported on the fixed set $S_1$. Then, following the same approach used in the case $p\in(2,4)$, we verify that \eqref{eq:Kest-p>4} implies  $\|K_{S_1}(\omega)\|_{L^p(d\omega)}\lesssim 1$.

Finally, we verify that the random set constructed in \eqref{eq:sp4} satisfies the uniform distribution as in \eqref{eq:slambda}. We observe that integer-valued independent random variables $\eta_{\zeta_j(\cdot)}(\cdot)$ are uniformly distributed over $\bar{I}_j$ for $1\leq j \leq L$, satisfying
\begin{align*}
    \P(\eta_{\zeta_j(\cdot)}(\cdot) =i)=N^{\frac{2}{p}-1} \quad \text{for} \quad i\in \bar{I_j}.
\end{align*}

\bigskip

\subsection{Preliminaries}

\begin{lem}[Lemma 4.4 in \cite{JuLaOrVu}]\label{te:probsum}
Let $A\subset [1,N]$, $|A|=l$ and $\max_{1\leq j \leq L}|A\cap I_j| \leq 1 $. Then we have
    \begin{align*}
        \|\sum_{i \in A} 1_{S(\omega)}(i)  \|_{L^q(d\omega)}\lesssim \frac{l}{s} + \frac{q}{\log(2+qs/l)}
    \end{align*}
    for any $q\geq 1$.
\end{lem}
\begin{proof}
Observe that $1_{S(\omega)}(i)$ for $i\in A$ are i.i.d. $\{0, 1\}$-valued random variables of mean $1/s$. We substitute $\xi_i$ and $\delta$ with $1_{S(\omega)}(i)$, $i\in A$ and $1/s$ respectively, and follow the proof of Lemma 4.4 in \cite{JuLaOrVu}.
\end{proof}
For $ \mathcal{E}\subset \R^N$ and $t>0$, let us denote the $L^2$ entropy of the set $\mathcal{E}$ at scale $t$ by $N_2(\mathcal{E},t)$. The entropy is defined by the minimal number of $L^2$ balls of radius $t$ to cover $\mathcal{E}$, that is,
\begin{align*}
    N_2(\mathcal{E},t) := \left\{N\ge 1: \exists \mathbf x_j \in \R^n, 1\le j \le n, \mathcal{E}\subset \bigcup_{j=1}^n(\mathbf x_j+tB_X)\right\}
\end{align*}
where $B_X \subset \R^N$ is the closed unit ball in $L^2$.

\begin{lem} [Lemma 4.1 in \cite{JuLaOrVu}]\label{te:probentropy}
    Let $\mathcal{E}\subset \R^N_+$, $B=\sup_{x\in \mathcal{E}}|x|_{l^2_N}$, and $S(\omega)$ be as in \eqref{eq:slambda}. Then for $m\in [1,N]$, and $1\leq q_0 < \infty$, we have
    \begin{align*}
        \| \sup_{\substack{x\in \mathcal{E},  |A|\leq m\\ \max_{1\leq j\leq L}|A\cap I_j|\leq 1}}\sum_{i\in  A} 1_{S(\omega)}(i)x_i \|_{L^{q_0}(d\omega)} \lesssim  \left ( \frac{m}{s}+\frac{q_0}{\log s} \right)^{1/2} B + (\log s)^{-1/2}\int_0^B\sqrt{\log N_2(\mathcal{E},t)}dt,
    \end{align*}
    where $N_2$ refers to the $L^2$ entropy.
\end{lem}

\begin{proof}
Since $1_{S(\omega)}(i)$ for $i\in A$ are i.i.d. $\{0, 1\}$-valued random selectors with mean $1/s$, we apply Lemma \ref{te:probsum} above in place of Lemma 4.4 from \cite{JuLaOrVu}, and follow the proof of Lemma 4.1 in \cite{JuLaOrVu}.
\end{proof}

\begin{te}[$2< p < 4$]\label{te:k2p4}
There exists $\sigma>0$ such that for any $\omega_2, \omega_3$
\begin{align}\label{eqn:3.29-Bourgain}
	\|K_{m_1,m_2,m_3}(\omega_1,\omega_2,\omega_3)\|_{L^{q_0}(d\omega_1)}\lesssim (\delta m_3^{p/2-1}+\frac{m_2+m_3}{m_1})^{1/2}(1+K(\omega_2)+K(\omega_3))^{p-\sigma},
\end{align}
where $q_0=\log N \sim \log \frac{1}{\delta}, N^{2/p}=\delta N$, with $K_{m_1,m_2,m_3}(\omega_1,\omega_2,\omega_3)$ defined by \eqref{eq:K_{m_1,m_2,m_3}} and $K(\omega)=K_{S(\omega)}$ given by \eqref{eq:KS}.
\end{te}

\begin{proof}
Fix $\omega_2$ and $\omega_3$.
Let $\mathcal{E}=\big\{\big(|\langle \varphi_i, 	f_{\bm b,\omega_2}(1+|f_{\bm c,\omega_3}|)^{p-2}\rangle	|\big)_{i=1}^{n}:\bm b\in \Pi_{m_2}, \bm c\in \Pi_{m_3}\big\}$ where $\Pi_{m_1}$ and $\Pi_{m_2}$ are defined in \eqref{eq:Pim}. We can write
$$
K_{m_1,m_2,m_3}(\omega_1,\omega_2,\omega_3)= \frac{1}{\sqrt{m_1}}\sup_{\substack{|A|\le m_1\\\max_{1\leq j\leq L}|A\cap I_j|\leq  1}}\sup_{\mathbf{x}\in \mathcal{E}}\Big(\sum_{i\in A}1_{S(\omega_1)}(i)x_i\Big).
$$
Therefore, by substituting Lemma 4.1 from \cite{JuLaOrVu} with Lemma \ref{te:probentropy} above (applied with $m=m_1$), we obtain 
\begin{align*}
	\|K_{m_1,m_2,m_3}(\omega_1,\omega_2,\omega_3)\|_{L^{q_0}(d\omega_1)}&\lesssim \left[\delta^{1/2}+m_{1}^{-1/2}\right]B+m_{1}^{-\frac{1}{2}}(\log N)^{-\frac12}\int_{0}^B[\log N_2(\mathcal{E},t)]^{\frac{1}{2}}\,dt,
\end{align*}
where $B = \sup_{\mathbf{x}\in \mathcal E}|\mathbf{x}|$. Then, we proceed by applying the proof of Theorem 6.1 in \cite{JuLaOrVu}.
\end{proof}

We recall our setup for the case $p>4$. Let $p_1$ satisfy $p/2< p_1 < p$, and assume that the $\Lambda(p_1)$-property holds for the random subset $S_1=\{\eta_k(\bar{\omega})\}_{k=1}^{L_1}$.

Let $\delta'=1/s'=N^{2/p-2/p_1}$, $L=N^{2/p}$. For each $1\leq j \leq L$, let $\{k(j)\}= \{ k\in [1,L'] : I_{k(j)} \subset \bar{I_j}\}$ be a subset of indices with $|\{k(j)\}|=s'$. Let $\zeta_j: \Omega' \rightarrow \{k(j)\}$ be independent random variables on some probability space $\Omega'$, uniformly distributed by
\begin{align*}
\P(\zeta_j=k(j))=\frac{1}{s'}.  
\end{align*}
Let $\{I_j'\}_{j=1}^{L}$ be a partition of $S_1$,  where $I_j'=\{\eta_{k(j)}(\bar{\omega})\}_{k(j)}$ and $|I_j'|=s'$.

Now, define $S_\omega\subset S_1$ to be a random subset obtained using the random variables $\{\zeta_j\}_{j=1}^L$ by
\begin{align*}
    S_\omega=\{\eta_{\zeta_j(\omega)} (\bar{\omega})  :  1\leq j \leq L \}
\end{align*}
and define
$K_{S_1}(\omega) = \sup_{|\mathbf a|_{l^2}\le 1}\|\sum_{i\in S_\omega}a_i\varphi_i\|_{p}.$

\begin{te}[$p>4$]\label{te:kp4}
Let $p>4$ and $p_1$ satisfy $p/2< p_1 < p$. Assume that \eqref{eq:lambda-p1} holds for a subset $S_1\subset[1,N]$ of cardinality $|S_1|\sim N^{2/p_1}=L_1$.  Then, there exists $\sigma>0$ such that for any $\omega_2$, $\omega_3$
the estimate \eqref{eq:Kest-p>4} holds
\begin{equation*}
\|K^{S_1}_{m_1,m_2,m_3}(\omega_1,\omega_2,\omega_3)\|_{L^{q_0}(d\omega_1)}\lesssim (\delta' m_3^{p/p_1-1}+\frac{m_2+m_3}{m_1})^{1/2}(1+K_{S_1}(\omega_2)+K_{S_1}(\omega_3))^{p-\sigma},
\end{equation*}
with $\delta'=N^{2/p-2/p_1}$, where 
$$
K^{S_1}_{m_1,m_2,m_3}(\omega_1,\omega_2,\omega_3):=\sup_{\substack{A\subset S_1 \\ |A|\le m_1\\ \max_{1\leq j \leq L'} |A\cap I_j'|\leq 1}}\sup_{\bm b\in \Pi^{S_1}_{m_2}}\sup_{\bm c\in \Pi^{S_1}_{m_3}}\frac{1}{\sqrt{m_1}}\sum_{i\in A} \xi_i(\omega_1)|\langle \varphi_i, 	f_{\bm b,\omega_2}|f_{\bm c,\omega_3}|^{p-2}\rangle |
$$
and
\begin{align*}
&\Pi^{S_1}_{m}=\{\mathbf{a}=(a_i)_{i\in S_1}:|\mathbf{a}|\le 1, |\supp \mathbf{a}|\le m\ \quad \text{and}\quad \max_{1\leq j \leq L} |\supp \bm a \cap I_j'| \leq 1\}.		
\end{align*}
\end{te}
\begin{proof}
We substitute probabilistic Lemma 4.1 from \cite{JuLaOrVu} with Lemma \ref{te:probentropy} above. Then, we follow the proof of Theorem 6.2 from \cite{JuLaOrVu}.
\end{proof}

\end{document}